\documentclass[runningheads,a4paper]{llncs}

\usepackage{epsfig,amsmath,amssymb,subfigure,version,graphicx,fancybox,mathrsfs,amscd}
\usepackage{cases,bm,multirow}
\usepackage{algorithm,algorithmic}
\usepackage{epstopdf}
\usepackage{url,xcolor}
\usepackage[colorlinks,urlcolor=blue]{hyperref}

\usepackage{todonotes}
\usepackage{comment}

\graphicspath{{images/}}

\newtheorem{definit}{\bf Definition}

\newcommand{\Lips}{Lipschitz}

\title{Fast Iterative Algorithms for Blind Phase Retrieval: A survey}
\author{ Huibin Chang${}^1$\thanks{Corresponding author}
    \and Li Yang${}^1$ 
	\and Stefano Marchesini${}^2$
	}
\institute
{
 ${}^1$School of Mathematical Sciences, Tianjin Normal University, Tianjin, China\\ {\tt E-mail:changhuibin@gmail.com}\\
${}^{2}$SLAC National Laboratory, CA, USA\\ {\tt E-mail:smarchesini@gmail.com}
}

\date{September 2021}

\begin{document}

\maketitle
\abstract{In nanoscale imaging technique  and ultrafast laser, the reconstruction procedure is normally formulated as a blind phase retrieval (BPR) problem, where one has to recover both the sample and the probe (pupil) jointly from phaseless  data. This survey  first presents the mathematical formula of BPR, related nonlinear optimization problems and then  gives a brief review of the recent iterative algorithms. It mainly consists of three types of algorithms, including the operator-splitting based first-order optimization methods, second order algorithm with Hessian, and subspace methods.  The future research directions  for experimental issues and theoretical analysis are further discussed.   }

\section{Introduction}
Phase retrieval (PR)  plays a key role in nanoscale imaging technique  \cite{Pfeiffer2018xray,elser2018benchmark,zheng2021concept,doga2022lensless} and ultrafast laser \cite{trebino1997measuring}. Retrieving the images of the sample from phaseless data is a long-standing problem. Generally speaking, designing fast and reliable algorithms is  challenging since directly solving the quadratic polynomials of PR is NP hard and the involved optimization problem is nonconvex and possibly nonsmooth. Thus, it has drawn the attentions of researchers for several decades \cite{Luke2005,Shechtman2014,Grohs2020phase,fannjiang_strohmer_2020}.  Among the general PR problems, besides the recovery of the sample, it is also of great importance to reconstruct the probes. The 
motivation of blind recovery is two-fold: (1) characteristics of the probe (wave front sensing); (2) improving the reconstruction quality of the sample. 
Essentially in practice, as the probe is almost never completely known, one has to solve such blind phase retrieval (BPR) problem, e.g. in coherent diffractive imaging (CDI)  \cite{thibault2012maximum}, convention ptychography imaging \cite{thibault2009probe,maiden2009improved}, Fourier ptychography \cite{zheng2013wide,ou2014embedded}, convolutional phase retrieval \cite{ahmed2018blind}, frequency-resolved optical gating \cite{trebino1997measuring} and others.

An early work by Chapman \cite{chapman1996phase} to solve the blind problem used the  Wigner-Distribution deconvolution method to retrieve the probe. In the optics community,    alternating projection (AP) algorithms are very popular for nonblind PR problems~\cite{marchesini2007invited,elser2018benchmark}. Some AP algorithms have also been applied to  BPR problems, \emph{e.g.} Douglas-Rachford (DR) based algorithm \cite{thibault2009probe}, extended ptychographic engine (ePIE) and variants  \cite{maiden2009improved,maiden2017further}, and relaxed averaged alternating reflections \cite{Luke2005} based projection algorithm  \cite{marchesini2016sharp}. More advanced first-order optimization method including proximal algorithms, \cite{hesse2015proximal,yan2020ptychographic,huang2021ptychography},  alternating direction of multipliers methods (ADMMs) \cite{chang2018Blind,fannjiang2020fixed} and convex programming method \cite{ahmed2018blind}.
To further accelerating the first-order optimization, several second-order algorithms utilizing the Hessian have also been developed  \cite{qian2014efficient,yeh2015experimental,ma2018globally,gao2017phaseless,kandel2021efficient}.  Moreover, the subspace methods \cite{CSIAM-AM-2-585} were successfully applied to the BPR as \cite{thibault2012maximum,chang2018Blind,fung2020multigrid}. 
\vskip .05in

The purpose of the survey is  to give a brief review of the recent iterative algorithms for BPR problem, so as to provide instructions for practical use and draw attentions of applied mathematician for further improvement.  The reminder of the survey is organized as follows: Section \ref{sec2} gives the mathematical formula for BPR and related nonlinear optimization  models, as well as  the closed-form expression of the proximal mapping.
Fast iterative algorithms  are  reviewed in section \ref{sec3}. Section \ref{sec4} further discusses the experimental issues and theoretical analysis. Section~\ref{sec5} summarizes this survey.

\section{Mathematical formula and nonlinear optimization model for BPR}\label{sec2}

First introduce the general nonblind PR problem in the discrete setting. 
By introducing a linear operator $A\in \mathbb C^{m,n}$, for the sample of interest $u\in\mathbb C^m$, experimental instruments usually collect the quadratic phaseless data $f\in\mathbb R^m$ as below:
\begin{equation}
\label{eq:model-GPR}
f=|Au|^2,
\end{equation}
in the ideal situation. 
However, noise contamination is evitable in practice \cite{chang2016Total} as
\begin{equation}\label{eq:model-PoiPR}
f_{noise}=\mathrm{Poi}(|Au|^2),
\end{equation}
where $\mathrm{Poi}$ denotes the random variable following i.i.d  Poisson distribution. See more advanced models for practical noise as outliers and structured and randomly distributed uncorrelated noise sources in \cite{godard2012noise,reinhardt2017beamstop,wang2017background,odstrvcil2018iterative,chang2019advanced} and references therein.

\subsection{Mathematical formula}

State the BPR problem starting from convention ptychography \cite{rodenburg2008ptychography}, since the principle of other BPR problems can be explained in a similar manner, all of which can  be unified as the  blind recovery problem. 

\begin{figure}[htb!]
\begin{center}
\includegraphics[width=.5\textwidth]{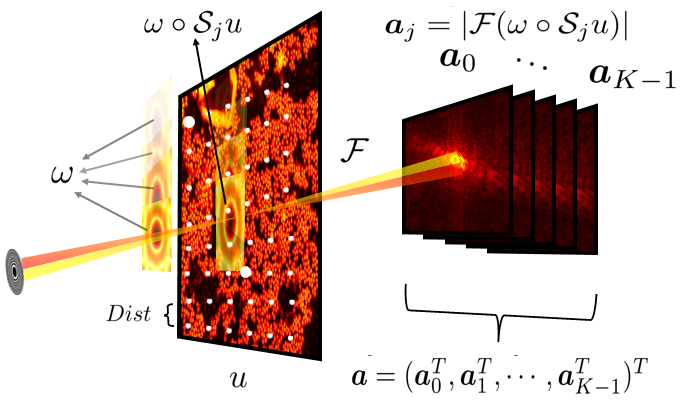}
\end{center}
\caption{Ptychographic phase retrieval (Far-field): A stack of phaseless data $f_j:=\bm a^2_j$ is collected, with
 $w$ being the localized coherent probe, and $u$ being the image of interest (sample).
 The white dots represent the scanning lattice points, with $Dist$ denoting the sliding distance between centers of two adjacent frames.}
\label{fig1}
\end{figure}

As shown in Figure \ref{fig1},  a detector in the far field measures  a series of phaseless intensities, by letting a localized coherent X-ray probe $w$  scan through the sample $u$. Let the 2D image and the localized 2D  probe   denote as $u\in\mathbb C^n$ with $\sqrt{n}\times\sqrt{n}$ pixels and  $w\in\mathbb C^{\bar m}$ with  $\sqrt{\bar m}\times \sqrt{\bar m}$ pixels, respectively. Here both the sample and the probe  are  rewritten as  vectors  by a lexicographical order. Let $f_j^P\in \mathbb R_+^{\bar m}~\forall 0\leq j\leq J-1$ denote the phaseless measurements satisfying 
\begin{equation}\label{eq:PtychoMeasurement}
f_j^P=|\mathcal F(w\circ \mathcal S_j u)|^2,
\end{equation}
where the symbols $|\cdot|$, $(\cdot)^2$, and $\circ$ and  represent the element-wise absolute value and square of a vector, and the element-wise multiplication of two vectors respectively, the symbol  $\mathcal S_j\in \mathbb R^{\bar m\times n}$ represents a matrix with binary elements extracting  a  patch  (with the index $j$ and  size  $\bar m$) from the entire sample, and the symbol $\mathcal F$ denotes the normalized discrete Fourier transformation (DFT). 
In practice, to get an accurate estimate of the probe, one has to solve a blind ptychographic PR problem.
Note that the coherent CDI problem \cite{thibault2012maximum} can be interpreted as a special blind ptychography problem with only one scanned  frame ($J=1$). 

A  recent super-resolution technique based on visible light called as Fourier ptychography method (FP) has been developed by Zheng et al. \cite{zheng2013wide} and quickly spreads out for fruitful applications \cite{zheng2021concept}. Letting $w$ and $u$ (Here reuse the notations for simplicity) being the  the point spread function (PSF) of the imaging system and the sample of interest, the collected phaseless data $f^{FP}_j$ of FP  can be expressed as 
\[
f^{FP}_j=|\mathcal F^{-1} (\bar w\circ \mathcal S_j \bar u)|^2~~\text{for~}0\leq j\leq J-1
\]
with $\bar w:=\mathcal F w$ and $\bar u:=\mathcal F u.$

Some similar problems dubbed as ``convolutional  PR'' were recently studied \cite{qu2017convolutional,qu2019convolutional,ahmed2018blind}. Given the sample $u$ and the convolution kernel $\kappa$, the phaseless measurement $f^{Cov}$ \cite{qu2017convolutional,qu2019convolutional} is given as
\[
f^{Cov}=|\kappa \circledast u|^2
\]
or 
\begin{equation}\label{Cov-1}
f^{Cov}=|\mathcal F(\kappa \circledast u)|^2, 
\end{equation}
where the symbol $\circledast$ denotes the convolution.

Other interesting  blind problem for full characterization of ultrashort optical 
pulses  is to use frequency-resolved optical gating (FROG) \cite{trebino1997measuring,bendory2017uniqueness,kane2021review}. The phaseless measurement for a typical SHG-FROG can be obtained as
\[
f^{FROG}_j=|\mathcal F(u\circ \mathcal T_j u)|^2,
\]
where the symbol $\mathcal T_j $ denotes the translation. From the measurement $\{\bm a_j^{FROG}\}_j$, one may also formulate it by BPR if assuming the element-wise multiplication for two independent variables.

All the mentioned problems  can be unified as the BPR problem, i.e. to recover the probe (pupil, convolution kernel or the signal itself) and the sample jointly. Essentially the relation between these two variables are bilinear. 
For conventional ptychography, the  bilinear operators  $\mathcal A:\mathbb C^{\bar m}\times \mathbb C^{n}\rightarrow \mathbb C^{m}$ and
 $\mathcal A_j:\mathbb C^{\bar m}\times \mathbb C^{n}\rightarrow \mathbb C^{\bar m}~\forall 0\leq j\leq J-1$,  are denoted as follows:
  \begin{equation}\label{eq:AAj}
 \mathcal A(w,u):=(\mathcal A_0^T (w,u), \mathcal A_1^T(w,u),\cdots, \mathcal A_{J-1}^T(w,u))^T,
 \end{equation}
 with 
 \[\mathcal A_j(w,u):=\mathcal F(w\circ \mathcal S_j u)\]
 and \[
 f:=(f_0^T, f_1^T, \cdots,  f_{J-1}^T)^T\in \mathbb R^m_+.\]
 Actually for all BPR problems, the bilinear operators can be unified as
 \begin{equation}\label{eq:bilinear}
\mathcal A_j(w,u):=\left\{
\begin{aligned}
&\mathcal F(w\circ \mathcal S_j u);&&\text{Case I: CDI and ptychography}\\
&\mathcal F^{-1}(\mathcal Fw\circ \mathcal S_j(\mathcal Fu));&&\text{Case II: Fourier ptychography}\\
&\mathcal F(w\circ \mathcal T_j u);&&\text{Case III: FROG}\\
&w\circledast u,\text{\ \ \ or\ \ \ }\mathcal F(w\circledast u);&&\text{Case IV: Convolution PR}\\
\end{aligned}
     \right.
 \end{equation}
 where there are totally one frame as $J=1$ for the last case for convolution PR.
 Hence by introducing the general bilinear operator 
 $\mathcal A(\cdot,\cdot)$, the BPR can be given below:
 \begin{equation}\label{equations:blindPR}
\text{BPR:\ }\text{To find the ``probe''~} w \text{~and the sample~} u,  ~s.t. ~ |\mathcal A(w, u)|^2=f, 
 \end{equation}
where $\mathcal A$ is denoted as \eqref{eq:AAj} and \eqref{eq:bilinear},
 and the per frame of phaseless measurements $f_j$ represents the measurement from four  cases.
Note that  the BPR problem is not limited to the cases with forward propagation as  \eqref{eq:bilinear}.

Denote two linear operators $A_w, A_u$ as below:
\begin{equation}\label{eq:linearOperator}
\begin{aligned}
&A_w u=\mathcal A(w,u)\forall u;\\
&A_u w=\mathcal A(w,u)\forall w;\\
\end{aligned}
\end{equation}
Then one can obtain the conjugate operators  
\begin{equation}\label{eq:AwT}
A_w^*z=\left\{
\begin{aligned}
&\sum_j \mathcal S_j^T(\mathrm{conj}(w)\circ \mathcal F^{-1}z_j);&&\text{Case I}\\
& \mathcal F^{-1}\sum_j \mathcal S_j^T(\mathrm{conj}(\mathcal F w)\circ \mathcal Fz_j);&&\text{Case II}\\
&\sum_j \mathcal T_j^T(\mathrm{conj}(w)\circ \mathcal F^{-1}z_j);&&\text{Case III}\\
&\mathrm{conj}(w)\circledast z,\text{\ \ \ or\ \ \ }\mathrm{conj}(w)\circledast \mathcal F^{-1}z;&&\text{Case IV}\\
\end{aligned}
\right.
\end{equation} 
and
\begin{equation}\label{eq:AuT}
A_u^*z=\left\{
\begin{aligned}
&\sum_j (\mathrm{conj}(\mathcal S_j u)\circ \mathcal F^{-1}z_j);&&\text{Case I}\\
& \mathcal F^{-1}\sum_j(\mathrm{conj}( \mathcal S_j\mathcal F u)\circ \mathcal Fz_j);&&\text{Case II}\\
&\sum_j (\mathrm{conj}(\mathcal T_j u)\circ \mathcal F^{-1}z_j);&&\text{Case III}\\
&\mathrm{conj}(u)\circledast z,\text{\ \ \ or\ \ \ }\mathrm{conj}(u)\circledast \mathcal F^{-1}z;&&\text{Case IV}\\
\end{aligned}
\right.
\end{equation} 
$\forall z=(z_1^T,z_2^T,\cdots,z_{J-1}^T)^T\in\mathbb C^m.$ Here $\sum_j$ is a simplified form of $\sum\nolimits_{j=0}^{J-1}.$
Consequently,  one obtains
\begin{equation}\label{eq:AtAw}
A_w^*A_wu=\left\{
\begin{aligned}
&\big(\sum_j \mathcal S_j^T |w|^2\big)\circ u;&&\text{Case I}\\
& \mathcal F^{-1}\big(\big(\sum_j \mathcal S_j^T |\mathcal F w|^2\big)\circ \mathcal F u\big);&&\text{Case II}\\
&\big(\sum_j \mathcal T_j^T |w|^2\big)\circ u;&&\text{Case III}\\
&\mathrm{conj}(w)\circledast w \circledast u;&&\text{Case IV}\\
\end{aligned}
\right.
\end{equation} 
and
\begin{equation}\label{eq:AtAu}
A_u^*A_u w=\left\{
\begin{aligned}
&\big(\sum_j \mathcal S_j |u|^2\big)\circ w;&&\text{Case I}\\
& \mathcal F^{-1}\big(\big(\sum_j \mathcal S_j |\mathcal F u|^2\big)\circ \mathcal F w\big);&&\text{Case II}\\
&\big(\sum_j \mathcal T_j |u|^2\big)\circ w;&&\text{Case III}\\
&\mathrm{conj}(u)\circledast u \circledast w.&&\text{Case IV}\\
\end{aligned}
\right.
\end{equation}

\subsection{Optimization problems and proximal mapping }

Solving a nonblind problem may be NP hard if knowing $w$ or $u$ in advance. 
Other than directly solving equations as  \eqref{equations:blindPR}, one can solve the following nonlinear optimization problems in order to determine the underlying image $u$ and probe $w$ from noisy measurements $f$:
\begin{equation}\label{eqPB}
\min\limits_{w,u}\mathcal M(|\mathcal A(w,u)|^2,f),
\end{equation}
where the symbol $\mathcal M(\cdot,\cdot)$ represents the error between the unknown intensity $|\mathcal A(w,u)|^2$  and collected phaseless data $f$. Various metrics  proposed under different noise settings include  amplitude based metric for Gaussian measurements (AGM) \cite{wen2012,chang2016phase}, intensity based metric for Poisson measurements (IPM) \cite{thibault2012maximum,chen2015solving,chang2016Total},  and intensity based metric for Gaussian measurements (IGM) \cite{qian2014efficient,candes2015phaseIT,sun2016geometric}, all of which 
 can be  expressed as
\begin{equation}
\mathcal M(g,f):=\left\{
\begin{aligned}
&\frac{1}{2}\big\|\sqrt{g}-\sqrt{f}\big\|^2;          &&\mbox{(AGM)} \\
&\frac{1}{2}\langle g-f\circ\log(g),\mathbf{1}\rangle;&&\mbox{(IPM)}\\
&\frac{1}{2}\|g-f\|^2;                                &&\mbox{(IGM)}\\
\end{aligned}
\right.
\label{eqDF}
\end{equation}
where the operations on vectors such as $\sqrt{\cdot}, \log(\cdot), |\cdot|, (\cdot)^2$ are all defined pointwisely in this survey,  $\mathbf{1}$ denotes a vector whose entries  all equals to ones, and $\|\cdot\|$ denotes the $\ell^2$ norm in Euclidean  space.

The proximal mapping for functions defined on complex Euclidean space is introduced below.
\begin{definit}\label{def3}
Given function $h: \mathbb C^N\rightarrow \mathbb R\bigcup \{+\infty\}$, the proximal mapping $\mathrm{Prox}_{h;\mu}:\mathbb C^N\rightarrow \mathbb C^N$ of $h$ is defined by
\begin{equation}\label{eq:prox-def}
\mathrm{Prox}_{h;\beta} (v) = \arg\min_{x} \Big(h(x)+\frac {\beta} 2 \|x-v\|^2\Big),
\end{equation}
with the symbol $\|\cdot\|$ denoted as
the $\ell^2$ norm of a complex vector  on complex Euclidean space (use the same notation for real and complex spaces). 
\end{definit}
Namely, the proximal operator for the function $\mathcal M(|\cdot|^2,f)$ defined in \eqref{eqDF} has a closed-form formula  \cite{chang2018variational} as below:
\begin{equation}
\mathrm{Prox}_{\mathcal M(|\cdot|^2,f);\beta}(z)=\left\{
\begin{aligned}
&\dfrac{\sqrt{f}+\beta|z|}{1+\beta}\circ \mathrm{sign}(z),                            &&\mbox{for AGM}; \\
&\dfrac{\beta |z|+\sqrt{(\beta|z|)^2+4(1+\beta)f}}{2(1+\beta)}\circ\mathrm{sign}(z), &&\mbox{for IPM} ; \\
&\varpi_{\beta}(|z|)\circ \mathrm{sign}(z),                                            &&\mbox{for IGM} ;
\end{aligned}
\right.
\label{eqProxCF}
\end{equation}
where $\forall~z\in \mathbb C^m,$  $(\mathrm{sign}(z))(t):=\mathrm{sign}(z(t))~\forall 0\leq t\leq m-1$,
$\mathrm{sign}(x)$ for a scalar $x\in \mathbb C$ is denoted as  $\mathrm{sign}(x)=\tfrac{x}{|x|}$ if $x\not=0,$ otherwise $\mathrm{sign}(0):=c$ with an arbitrary constant $c\in \mathbb C$ with unity length, and 
\begin{equation}\label{cubicsolver}
\begin{split}
\varpi_\beta(|z|)(t)=\left\{
\begin{split}
&\sqrt[3]{\tfrac{\beta|z(t)|}{4}+\sqrt{D(t)}}+\sqrt[3]{\tfrac{\beta|z(t)|}{4}-\sqrt{D(t)}},~\mbox{if~} D(t)\geq 0; \\
&2\sqrt{\tfrac{f(t)-\tfrac{\beta}{2}}{3}}\cos\Big(\arccos\tfrac{\theta(t)}{3}\Big),\qquad \mbox{otherwise,}
\end{split}
\right.
\end{split}
\end{equation}
 for  $0\leq t\leq m-1,$
 with $
D(t)=\frac{(\tfrac{\beta}{2}-f(t))^3}{27}+\dfrac{\beta^2|z(t)|^2}{16},$
and
$\theta(t)=\dfrac{\beta|z(t)|}{4\sqrt{\frac{(f(t)-\frac{\beta}{2})^3}{27}}}.$

Note that  the alternating direction method of multipliers (ADMM) was adopted in \cite{wen2012,chang2016phase,chang2016Total} to solve the variational phase retrieval model in \eqref{eqPB}. However, due to the lack of the globally  \Lips~ differentiable terms in the objective function,  it seems difficult to guarantee its convergence. Some other variants of the metric have been recently proposed, such as the penalized metrics $\mathcal M(|\cdot|^2+\epsilon\bm 1, f+\epsilon\bm 1)$ by adding a small positive scalar $\epsilon$ as \cite{guizar2008phase,chang2018Blind,gao2019solving}. 
Although it has simple form,  the   technique will make the related proximal mapping  not have closed form expression,  such that  additional computation cost as an inner loop may have to be introduced \cite{chang2018Blind}. 
By cutting off the AGM near the origin, and then adding back a smooth function, one can  keep the global minimizer unchanged. Hence, a novel smooth truncated AGM (ST-AGM) $\mathcal G_\epsilon(\cdot; f)$ with truncation parameter $\epsilon>0$ \cite{chang2021over} was  designed below:
\begin{equation}\label{eq:STAGM}
\mathcal M_\epsilon(z,f):=\sum\nolimits_j M_\epsilon(z(j),f(j)),
\end{equation}
where $\forall\ \  x\in \mathbb C, b\in\mathbb R^+,$
\begin{equation}\label{eq:obj-1}
M_{\epsilon}(x, b):=
\left\{
\begin{aligned}
&\frac{1-\epsilon}{2}\left( b-\tfrac{1}{\epsilon}{|x|^2}\right), \text{\ \ if\ } |x|<\epsilon \sqrt{b};\\
&\frac{1}{2} \big||x|-\sqrt{b}\big|^2, \text{\ \ \ \ \ otherwise.}
\end{aligned}
\right.
\end{equation}
Readily its closed form of the corresponding proximal mapping is given directly as
\begin{equation}
\begin{aligned}
&(\mathrm{Prox}_{\mathcal M_\epsilon;\beta}(y))(j)\\
=&\mathrm{sign}(y(j))\times
\left\{
\begin{aligned}
&\max\Big\{0, \tfrac{\beta|y(j)|}{\beta-\tfrac{1-\epsilon}{\epsilon}}\Big\},&\text{if~} |y(j)|<\epsilon-\tfrac{1-\epsilon}{\beta}\sqrt{f(j)};\\
&\tfrac{\sqrt{f(j)}+\beta|y(j)|}{1+\beta},&\text{otherwise,}\\
\end{aligned}
\right. 
\end{aligned}
\label{eq:Prox_sol}
\end{equation}
if $\epsilon\in (0,1)$.  More other elaborate metrics can be found \cite{Luke2005,cai2021solving} and references therein.

\section{Fast iterative algorithms }\label{sec3}

In this section,   the main iterative algorithms for BPR will be introduced. Note that each algorithm may be designed originally for a specific case of \eqref{eq:bilinear}. Hence   the basic idea based on the original case will be explained first, and  the possible extensions to other cases will be discussed then.

\subsection{Alternating projection (AP) algorithms }
First consider BPR defined in \eqref{equations:blindPR} in the case of convention ptychography. 

Given the  exit wave in the far-field $\Psi:=(\Psi^T_0,\Psi^T_1,\cdots,\Psi_{J-1}^T)^T\in \mathbb C^{m},$ with 
\[\Psi_j:=\mathcal F(w\circ\mathcal S_j u)~\forall 0\leq j\leq J-1,\]
the optimal  exit wave $\Psi^\star$ lies in the intersection of two following sets, \emph{i.e.}
\[
\Psi^\star\in\widehat{\mathscr X}_1\bigcap \widehat{\mathscr X}_2,
\]
with
\begin{equation}\label{eqX}
\begin{aligned}
&\widehat{\mathscr X}_1:=\{\Psi:=(\Psi^T_0,\Psi^T_1,\cdots,\Psi_{J-1}^T)^T\in \mathbb C^m:~|\Psi_j|=\sqrt{f_j}~\forall 0\leq j\leq J-1\},\\
&\widehat{\mathscr X}_2:=\{\Psi\in \mathbb C^m:~\exists w\in \mathbb C^{\bar m}, u\in \mathbb C^n, s.t. ~~w\circ \mathcal S_j u=\mathcal F^{-1}\Psi_j~\forall 0\leq  j\leq J-1 \}.
\end{aligned}
\end{equation}

The AP algorithm determining this intersection alternatively calculates the projections onto these two sets $\widehat{\mathscr X}_1$ and $\widehat{\mathscr X}_2$.
Regarding the projection onto $\widehat{\mathscr X}_1$ as
\[
\widehat {\mathcal P}_1(\Psi):=\arg\min_{\widehat\Psi\in \widehat {\mathscr X}_1}\|\widehat \Psi-\Psi\|^2,
\]
 one readily gets  a closed form solution
\[\widehat {\mathcal P}_1(\Psi):=\big((\widehat {\mathcal P}_1^0(\Psi))^T,\cdots, (\widehat {\mathcal P}_1^{J-1}(\Psi))^T\big)^T,
\]
with
\[
\widehat {\mathcal P}_1^j(\Psi)=\sqrt{f_j}\circ \mathrm{sign}(\Psi_j)~0\leq j\leq J-1.
\]
For the projection onto $\widehat{\mathscr X}_2$, given $\Psi^k$ as the solution
in the $k^{\text{th}}$ iteration, one gets
\[
\begin{aligned}
\widehat {\mathcal P}_2(\Psi^{k}):=((\mathcal F(w^{k+1}\circ \mathcal S_0 u^{k+1}))^T, &(\mathcal F(w^{k+1}\circ \mathcal S_1 u^{k+1}))^T,\cdots,\\
&\qquad (\mathcal F(w^{k+1}\circ \mathcal S_{J-1} u^{k+1}))^T)^T,
\end{aligned}
\]
where
\begin{equation}\label{modelLS}
(w^{k+1}, u^{k+1})=\arg\min\limits_{w, u} F(w,u,\Psi^k):=\tfrac12  {\textstyle \sum\nolimits_j}\|\mathcal F^{-1}\Psi^k_j-w\circ \mathcal S_j u\|^2.
\end{equation}
Unfortunately, it does not have a closed form solution.
One can  solve \eqref{modelLS}  by alternating minimization (with $T$ steps) as below:
\begin{equation}
\begin{split}
&w_{l+1}=\arg\min_{w} F(w, u_l,\Psi^k),  \\
&u_{l+1}=\arg\min_u F(w_{l+1},u,\Psi^k) ~\forall l=0,1,\cdots,T-1.
\end{split}
\label{eqAPmu}
\end{equation}
Readily one has
\begin{equation}
\begin{split}
&w_{l+1}\approx
\frac{\sum\nolimits_j\mathrm{conj}(\mathcal S_j u_l)\circ \mathcal F^{-1}\Psi_j^k}
{\sum\nolimits_j\left|\mathcal S_j u_l\right|^2+\bar\alpha_1};\\
&u_{l+1}\approx\frac{\sum_j \mathcal S_j^T (\mathrm{conj}(w_{l+1})\circ \mathcal F^{-1}\Psi_j^k)}{ \sum_j (\mathcal S_j^T |w_{l+1}|^2) +\bar\alpha_2} ~\forall l=0,1,\cdots,T-1,
\end{split}
\label{eq:subsolverAP}
\end{equation}
where the parameters $0<\bar\alpha_1, \bar \alpha_2\ll 1$ are introduced in order to avoid dividing by zeros. 

Letting $\Psi^k$ be iterative solution in the $k^{th}$ iteration, the standard AP  for BPR can be directly given as below:

\vspace{2mm}
\hspace{1mm}(1) Compute $\widehat\Psi^k$ by $\widehat\Psi^k_j=\mathcal F(w^{k+1}\circ\mathcal S_j u^{k+1})$, where the pair $(w^{k+1},u^{k+1})$ is approximately  solved  by \eqref{eq:subsolverAP}.

\hspace{1mm}(2) Compute $\Psi^{k+1}$ by
\vspace{3mm}
$
\Psi^{k+1}=\widehat {\mathcal P}_1(\Psi^k).
$

The DR algorithm for BPR can be formulated in two steps  \cite{thibault2009probe}, as follows:
\vspace{2mm}

\hspace{1mm}(1) Compute $\widehat\Psi^k$ as the first step of AP.

\hspace{1mm}(2) Compute $\Psi^{k+1}$ by
\vspace{3mm}
\begin{equation}\label{DR}
\Psi^{k+1}=\Psi^k+\widehat {\mathcal P}_1(2\widehat\Psi^k-\Psi^k)-\widehat\Psi^k.
\end{equation}
Note that the formula \eqref{DR} utilizing Douglas-Rachford operator is essentially Fienup’s hybrid input–output map,  which can also be derived with proper parameters from 
difference map \cite{Elser2003}.

Since the fixed point of DR iteration may not exist, Marchesini et al. \cite{marchesini2016sharp} adopted the relaxed version of DR (dubbed as RAAR by Luke \cite{Luke2005}) to further improve the stability of the reconstruction from noisy measurements, which simply takes weighted average of right term of \eqref{DR} and $\widehat\Psi^k$  with a tunable parameter $\delta\in (0, 1)$ as
\[
\Psi^{k+1}=\delta\big(\Psi^k+\widehat {\mathcal P}_1(2\widehat\Psi^k-\Psi^k)-\widehat\Psi^k\big)+(1-\delta)\widehat\Psi^k,
\]
with $\widehat\Psi^k$ determined in a same manner as the first step of AP.

At the end of this part, extension of AP to general BPR problems will be discussed. Similarly as for the ptychography, introduce $\Psi$
as
\[
\Psi=\mathcal A(w,u), \text{and~} 
\Psi_j=\mathcal A_j(w,u).
\] 
In a same manner, one can define two constraint sets and establish the AP algorithms for the four cases of BPR. The only differences lie in the calculations of the projections on to the bilinear constraint set. As  \eqref{modelLS}, consider
\begin{equation}\label{eq:LS-bilinear}
\min_{w,u}\|\Psi^k-\mathcal A(w,u)\|
\end{equation}
by alternating minimization, where $\Psi^k$ is the iterative solution. 
Then the scheme is given below:
\begin{equation}
\begin{split}
&w_{l+1}\approx (A_{u_l}^*A_{u_l}+\bar\alpha_1\mathbf I)^{-1}A_{u_l}^*\Psi^k,\\
&u_{l+1}\approx(A_{w_{l+1}}^*A_{w_{l+1}}+\bar\alpha_2\mathbf I)^{-1}A_{w_{l+1}}^*\Psi^k ~\forall l=0,1,\cdots,T-1.
\end{split}
\label{eq:subsolverAP-general}
\end{equation}
The detailed forms of these operators can be found in \eqref{eq:AwT}-\eqref{eq:AtAu}. Notably the inverse in \eqref{eq:subsolverAP-general} can be efficiently solved by point-wise division or DFT.

\subsection{ePIE-type algorithms }

This iterative algorithm can be expressed as an AP method for convention ptychography as follows: To find $\Psi^\star_{n_k}$ belonging to the intersection as
\[
\Psi^\star_{n_k}\in\{|\Psi_{n_k}|=\sqrt{f_{n_k}}\}\cap\{\Psi_{n_k}:~\exists w\in\mathbb C^{\bar m}, u\in \mathbb C^n, ~~s.t. ~~w\circ\mathcal S_{n_k} u=\mathcal F^{-1}\Psi_{n_k}\},
\]
with  a random frame index $n_k$.
Let $w^k, u^k$ be the iterative solutions in the $k^{\text{th}}$ iteration.
 By first computing the projection of $\psi^k_{n_k}:=\mathcal F(w^k\circ\mathcal S_{n_k}u^k)$  by $\widehat {\mathcal P}_1^{n_k}(\psi^k_{n_k})$, and then updating $w^{k+1}$ and $u^{k+1}$ by the gradient descent algorithm (inexact projection) for \eqref{modelLS},
the ePIE algorithm proposed by Maiden and Rodenburg \cite{maiden2009improved} can be expressed by
updating $w^{k+1}$ and $u^{k+1}$  in parallel as
\begin{equation}
\label{eq:ePIE}
\left\{
\begin{aligned}
&w^{k+1}=w^k-\frac{d_2}{\|\mathcal  S_{n_k} u^{k}\|^2_{\infty}}\mathcal S_{n_k} \mathrm{conj}(u^{k})\circ\mathcal F^{-1}(\Psi^{k}_{n_k}-\widehat {\mathcal P}_1^{n_k} (\Psi^{k}_{n_k}))\\
&u^{k+1}=u^k-\frac{d_1}{{\|\mathcal S_{n_k}^Tw^k\|^2_{\infty}}} \mathcal S_{n_k}^T\left({\mathrm{conj}(w^k)}\circ\mathcal F^{-1}(\Psi^{k}_{n_k}-\widehat {\mathcal P}_1^{n_k} (\Psi^{k}_{n_k}))\right),
\end{aligned}
\right.
\end{equation}
with frame index $n_k\in \{0,1,\cdots,J-1\}$   generated randomly, and positive parameters $d_1$ and $d_2$ (default values are ones)
 and  $\|w\|_\infty:=\max_t |w(t)|$.

The regularized PIE (rPIE) was further proposed by Maiden, Johnson and Li \cite{maiden2017further} 
as 
\begin{equation}
\label{eq:rPIE}
\left\{
\begin{aligned}
w^{k+1}=w^k-&\frac{\bm 1}{\delta\|\mathcal  S_{n_k} u^{k}\|^2_{\infty}+(1-\delta)\mathcal  S_{n_k} |u^{k}|^2 }\\
&\circ\mathcal S_{n_k} \mathrm{conj}(u^{k})\circ\mathcal F^{-1}(\Psi^{k}_{n_k}-\widehat {\mathcal P}_1^{n_k} (\Psi^{k}_{n_k})),\\
u^{k+1}=u^k-&\frac{\bm 1}{\delta{\|\mathcal S_{n_k}^Tw^k\|^2_{\infty}}+(1-\delta)\mathcal S_{n_k}^T|w^k|^2} \\
&\circ\mathcal S_{n_k}^T\left({\mathrm{conj}(w^k)}\circ\mathcal F^{-1}(\Psi^{k}_{n_k}-\widehat {\mathcal P}_1^{n_k} (\Psi^{k}_{n_k}))\right),
\end{aligned}
\right.
\end{equation}
with the scalar constant $\delta \in (0, 1)$. It can be interpreted as a hybrid scheme for the stepsize of gradient descent, which takes the weighed average of the denominator of the ePIE scheme \eqref{eq:ePIE}  and 
first term in the denominator of  AP scheme  \eqref{eq:subsolverAP}. The rPIE algorithm was further accelerated by momentum \cite{maiden2017further}.

One can directly get the ePIE and rPIE schemes  for FP \cite{zheng2021concept} by replacing the variables $w$ and $u$ by $\mathcal F w$ and $\mathcal F u$. 
The ePIE type algorithms are very popular in optics community, since it is enough to implement the algorithm if one knows how to calculate the gradient of the objective functions, and the memory footprint is much smaller than more advanced AP algorithm including DR and RAAR. However, it tends to unstable when the data redundancy is insufficient (e.g. noisy data, big-step scan) as reported in \cite{chang2018Blind}. Moreover, the theoretical convergence is unknown and seems challenging due to the relation with nonsmooth objective functions. 

Note that if with totally $J=1$ frame as CDI, the differences between the ePIE (with $d_1=d_2=1$) and standard AP lie in the preconditioning matrices: AP utilizes the spatial weighted diagonal matrices $A_u^*A_u$ and $A_w^*A_w$, while ePIE utilizes the spatial-independent constant determined by the maximum of their diagonal matrices.

\subsection{Proximal algorithms}\label{sec1-1-2}

\subsubsection{Proximal heterogeneous block implicit-explicit (PHeBIE)}
For convention ptychography, consider   an optimization  problem (To get rid of introducing redundant notations in this survey,  slightly modify the constraint set of $\Psi$ in \cite{hesse2015proximal} as $\widehat{\mathscr X_1}$ and adjust the notation of the first term of the following model accordingly in order to present  an equivalent form)  \cite{hesse2015proximal}  as follows:
 \begin{equation}\label{eqLS}
 \min\limits_{w,u,\Psi} F(w,u,\Psi)+\mathbb I_{\widehat{\mathscr X}_1}(\Psi)+{\mathbb I_{\mathscr X_1}(w)+\mathbb I_{\mathscr X_2}(u),},
 \end{equation}
 with $F(w, u, \Psi)$ and $\widehat{\mathscr X}_1$  denoted in \eqref{modelLS} and \eqref{eqX}, respectively,
 and the indicator function $\mathbb I_{{\mathscr X}}$  denoted as
 \[
 \mathbb I_{\mathscr X} (\Psi):=
 \left\{
 \begin{aligned}
 &0,~\text{if~}\Psi\in \mathscr X,\\
 &+\infty,~\text{otherwise},
 \end{aligned}
 \right.
\]
where
 the amplitude constraints of the probe and  image  are incorporated (In \cite{hesse2015proximal}, the authors further considered the compact support condition of the probes and image), where 
 \begin{equation}\label{eq:constr}
 \begin{split}
 &\mathscr X_1:=\{w\in\mathbb C^{\bar m}:~\|w\|_\infty\leq C_{w}\};\\
 &\mathscr X_2:=\{u\in \mathbb C^n:~\|u\|_\infty\leq C_{u}\}
 \end{split}
 \end{equation} 
 with two positive constants $C_w, C_u$.
The projection operator 
onto  $\mathscr X_1$ is readily obtained as
\[
\mathrm{Proj}( w;C_ w):=\min\{C_ w, | w|\}\circ \mathrm{sign}( w)~\forall w,
\]
which is the close form expression for the minimizer to the problem
\[
\min_{\|\tilde w\|_\infty\leq C_w} \tfrac12\|\tilde w- w\|^2.
\]
Similarly one gets the projection onto $\mathscr X_2$ as 
$\mathrm{Proj}( u;C_ u)$.

Hesse et al.  \cite{hesse2015proximal} further adopted the proximal alternating linearized minimization (PALM) method \cite{bolte2014proximal}  for the BPR problem in the case of convention ptychography, such that the  proximal heterogeneous block implicit-explicit (PHeBIE) [See \cite{hesse2015proximal}, Algorithm 2.1] consists of  two steps with  three positive parameters $d_1, d_2$ and $\gamma$:
\vspace{1.5mm}

\hspace{1mm}(1)  Calculate $w^{k+1}, u^{k+1}$ sequentially 
 as
\begin{equation}
\label{eq:PGM-step1}
\left\{
\begin{aligned}
w^{k+1}&{=}\mathrm{Proj}\Big(w^k-\tfrac{1}{d_1^k}\sum\nolimits_j \mathcal S_j \mathrm{conj}(u^k)\circ\mathcal F^{-1}(w^k\circ\mathcal S_j u^k- \Psi_j^k); C_w\Big);\\
u^{k+1}&{=}\mathrm{Proj}\Big(u^k-\tfrac{1}{d_2^k}\sum\nolimits_j\mathcal S_j^T(\mathrm{conj}(w^{k+1}) \circ \mathcal F^{-1}(w^{k+1}\circ \mathcal S_j u^k -\Psi_j^k)), C_u\Big),
\end{aligned}
\right.
\end{equation}
with $d_1^k:=d_1 \big\|\sum_j|\mathcal S_j u^k|\big\|_{\infty}^2$,  $d_2^k:=d_2\big\|\sum_j \mathcal S_j^T|w^{k+1}|\big\|_{\infty}^2$.
\vspace{1.5mm}

\hspace{1mm}(2) Calculate $\Psi^{k+1}$ by
\[
\Psi^{k+1}=\widehat {\mathcal P}_1\big(\tfrac{1}{{1+\gamma}}({\widehat\Psi^{k+1}+\gamma \Psi^k})\big),
\]
with \[
\widehat\Psi^{k+1}_j:=\mathcal F(w^{k+1}\circ S_j u^{k+1}).
\]

\vspace{1.5mm}

Under  the assumption of  boundedness of iterative sequences, the convergence of PALM to  stationary points of \eqref{eqLS} was proved~\cite{hesse2015proximal}. To the knowledge, it is the first iterative algorithm with convergence guarantee for BPR problem.  The PHeBIE  has multiple similarities with ePIE. The main differences between them  are that for ePIE, only the gradient of $F$ \emph{w.r.t.} a randomly selected single frame is adopted to update $w$ and $u$ per outer loop as \eqref{eq:ePIE} , while for PHeBIE, each block of $w$ and $u$ can be updated in parallel by employing the gradient as \eqref{eq:PGM-step1} \emph{w.r.t.} all adjacent frames. Therefore,  PHeBIE is  more stable than ePIE numerically.

One readily knows that 
the convergence rate relies on the \Lips~constant of partial derivative of $F$. In order to get smaller constant, a directly way is to employ the derivative of a small block for unknowns. Hence, based on the partition of the sample and the probe, the parallel version of PHeBIE  was also provided \cite{hesse2015proximal} with convergence guarantee.

For more extensions to other cases of BPR,
one can introduce a generalized nonlinear optimization model:
\[
\min_{w, u, \Psi}  \|\Psi-\mathcal A(w, u)\|^2+\mathbb I_{\widehat{\mathscr X}_1}(\Psi)+\mathbb I_{\mathscr X_1}(w)+\mathbb I_{\mathscr X_2}(u),
\]
where one adopts the  same form  as \eqref{eq:LS-bilinear} for the first term. The detailed algorithms are omitted, since one only needs to update the gradient of first term following \eqref{eq:AwT}-\eqref{eq:AtAu}.

\subsubsection{Variant of Proximal algorithm }
Here introduce a general constraint set
for the bilinear relation as 
\begin{equation}\label{eq:bilinearSet}
X:=\{\Psi\in \mathbb C^m:~\exists w\in \mathbb C^{\bar m}, u\in \mathbb C^n, s.t.~~ \mathcal A(w,u)=\Psi\}.    
\end{equation}

Consider the optimization problem as
\begin{equation}
\min_{z} \mathbb I_{\widehat{\mathscr X}_1}(z)+\mathbb I_{ X}(z).    
\end{equation}
By replacing the indicator function by the metrics, and further combining the alternating  minimization with proximal algorithms, Han  \cite{yan2020ptychographic} derived a new proximal  algorithm for the convention ptychography problem. 
Specifically, the proposed algorithm with a generalized form for BPR has the following steps:
\begin{equation}\label{eq:proximal-ptycho}
\begin{split}
&\text{Step 1:\quad} z^{k+1}=\arg\min_z \mathcal M(|z|^2;f)+\tfrac{\beta}{2}\|z-\mathcal A(w^k,u^k)\|^2; \\
&\text{Step 2:\quad}w^{k+1}=\arg\min_w \|z^{k+1}-A_{u^{k}} w)\|^2.\\
&\text{Step 3:\quad}u^{k+1}=\arg\min_u \|z^{k+1}-A_{w^{k+1}} u)\|^2.\\
\end{split}
\end{equation}
Here the last two steps can be solved in a same manner as \eqref{eq:subsolverAP-general}.
The above algorithm has deep connections with the ADMM \cite{chang2018Blind}. If  removing the constraint of boundedness of two variables, and setting the penalization parameter to zero in \eqref{eqModel}, then by solving the constraint problem \eqref{eqI-1} by adding a penalization term $\|z-\mathcal A(w,u)\|^2$ without introducing the multiplier $\Lambda$, one can get exactly the same iterative scheme as \eqref{eq:proximal-ptycho}.  
Besides, it was  further improved by 
accelerated proximal gradient method in \cite{yan2020ptychographic}, and recently by stochastic gradient descent  \cite{huang2021ptychography} for FP.

\subsection{ADMM}

As a typical operator-splitting algorithm, ADMM is very flexible and successfully applied to inverse and imaging problems \cite{Wu&Tai2010,boyd2011distributed}, which is also adopted for classical and ptychographic PR   problems \cite{wen2012,chang2018Blind}.

Consider the  metrics using penalized-AGM (pAGM) and penalized-IPM (pIPM) as  to measure the error of recovered intensity and the targets. A nonlinear  optimization model \cite{chang2018Blind}  was given as:
\begin{equation}\label{eqModel}
 \min\nolimits_{w\in \mathbb C^{\bar m},u\in\mathbb C^n}\mathcal G(\mathcal A(w,u)){+\mathbb I_{\mathscr X_1}(w)+\mathbb I_{\mathscr X_2}(u),}
\end{equation}
with $\mathcal G(z):=\mathcal M(|z|^2+\epsilon \bm 1,f+\epsilon\bm 1)$ and the constraint sets defined in \eqref{eq:constr}.
The authors further
leveraged the additional data $\bm c\in\mathbb R^{\bar m}_+$ to eliminate structural artifacts caused by grid scan, and therefore obtained the following variant
\begin{equation}\label{eqModelpriorII}
\min\limits_{ w\in \mathbb C^{\bar m},u\in\mathbb C^n}
 \mathcal G(\mathcal A( w,u)){+\mathbb I_{\mathscr X_1}( w)+\mathbb I_{\mathscr X_2}(u)}+\tau \widehat{\mathcal G}(\mathcal F w),
\end{equation}
where $\tau$ is a positive parameter,
{the additional measurement $\bm c$ is the diffraction pattern (absolute value of Fourier transform of the probe) as $\bm c:=|\mathcal F u|$, and
$
\widehat{\mathcal G}(z):=\mathcal B(|z|^2+\epsilon\bm 1,\bm c^2+\epsilon\bm 1).
$
}
For simplicity, assume that $\mathcal G$ and $\hat{\mathcal G}$ adopt the same metric.

As the procedures for solving the two above models are quite similar using ADMM,  only details for solving  the first optimization model \eqref{eqModel} is given below. 
By introducing an auxiliary variable
$z=\mathcal A(w, u)\in\mathbb C^{m}$,   an  equivalent form of \eqref{eqModel} is formulated below:
\begin{equation}\label{eqI-1}
\min\limits_{w,u,z} \mathcal G(z){+\mathbb I_{\mathscr X_1}(w)+\mathbb I_{\mathscr X_2}(u),}~~s.t.~~~z-\mathcal A(w, u)=0.
\end{equation}
The corresponding augmented Lagrangian reads
\begin{equation}\label{eqAL}
\begin{aligned}
&\Upsilon_\beta (w,u,z,\Lambda):=\mathcal G(z){+\mathbb I_{\mathscr X_1}(w)+\mathbb I_{\mathscr X_2}(u)}
&+&\Re(\langle z-\mathcal A(w, u), \Lambda \rangle)
\\
&&+&\tfrac{\beta}{2}\|z-\mathcal A(w, u)\|^2,
\end{aligned}
\end{equation}
with the multiplier $\Lambda\in \mathbb C^m$ and  a positive parameter $\beta,$
where  $\Re(\cdot)$ denotes the real part of a complex number. Then one considers the following problem:
\begin{equation}\label{eqSaddle}
\max_{\Lambda}\min\limits_{w,u,z} \Upsilon_\beta(w,u,z,\Lambda).
\end{equation}
 Given the approximated solution $(w^k,u^k,z^k,\Lambda^k)$ in the $k^{\text{th}}$ iteration, the four-step iteration by the generalized ADMM (only the subproblems \emph{w.r.t.} $w$ or $u$ have proximal terms is given as follows:
\begin{equation}\label{eqADMM}
\left\{
\begin{aligned}
&\mbox{Step 1:\quad} w^{k+1}=\arg\min\limits_w {\Upsilon_\beta(w,u^k,z^k,\Lambda^k)+\tfrac{\alpha_1}{2}\|w-w^k\|^2_{M_1^k}}\\
&\mbox{Step 2:}\quad u^{k+1}=\arg\min\limits_u  {\Upsilon_\beta(w^{k+1},u,z^k,\Lambda^k)+\tfrac{\alpha_2}{2}\|u-u^k\|^2_{M_2^k}}\\
&\mbox{Step 3:\quad} z^{k+1}=\arg\min\limits_z \Upsilon_\beta(w^{k+1},u^{k+1},z,\Lambda^k)\\
&\mbox{Step 4:\quad} \Lambda^{k+1}=\Lambda^{k}+\beta(z^{k+1}-\mathcal A(w^{k+1},u^{k+1})),\\
\end{aligned}
\right.
\end{equation}
with  diagonal positive semidefinite  matrices  $M^k_1\in\mathbb R_+^{\bar m\times\bar m}$ and $M^k_2\in \mathbb R_+^{n\times n}$ and two penalization parameters $\alpha_1, \alpha_2>0,$
where $\|w\|^2_{M_1^k}:=\langle M_1^kw,w\rangle$ and
$\|u\|^2_{M_2^k}:=\langle M_2^ku,u\rangle$.

Detailed algorithms will be given focusing on convention ptychography as \cite{chang2018Blind}. Note that these two matrices $M_1^k, M_2^k$ are assumed to be diagonal so that subproblems in Step 1 and Step 2 have closed form solutions.
Roughly speaking, based on splitting technique of proximal ADMM, subproblems of $u$, $w$ and $z$ are element-wise optimization problems with closed form solutions, such that  each subproblem  can be fast solved.
In practice, these two matrices are chosen by hand,  and an adaptive strategy  was presented in \cite{chang2018Blind} in order to guarantee the convergence.
Letting 
\[
\hat z^k:=z^k+\tfrac{\Lambda^k}{\beta}.
\]
{and  the diagonal matrices $M_1^k$ and $M_2^k$  satisfy
\begin{equation}\label{eqM1-Cond}
\left\{
\begin{aligned}
&\min\nolimits_t  \sum\nolimits_j\left|\left(\mathcal S_j u^k\right)(t)\right|^2+\tfrac{\alpha_1}{\beta}\mathrm{diag}(M_1^k)(t)>0,\\
&\min\nolimits_t  \sum\nolimits_j\left|\left(\mathcal S_j^T  w^{k+1}\right)(t)\right|^2+\tfrac{\alpha_2}{\beta}\mathrm{diag}(M_2^k)(t)>0,
\end{aligned}
\right.
\end{equation}
the close form solutions
 of Step 1 and Step 2 are given  as
\begin{equation}\label{eqOmega}
\left\{
\begin{aligned}
& w^{k+1}=\mathrm{Proj}\Big(
\tfrac{\beta\sum\nolimits_j\mathrm{conj}(\mathcal S_j u^k)\circ(\mathcal F^{-1} \hat z_j^k){+{\alpha_1} \mathrm{diag}(M_1^k)\circ w^k}}
{\beta\sum\nolimits_j\left|\mathcal S_j u^k\right|^2+{\alpha_1}\mathrm{diag}(M_1^k)}
; C_w\Big);\\
&u^{k+1}=\mathrm{Proj}\Big(\tfrac{\beta\sum_j \mathcal S_j^T (\mathrm{conj}( w^{k+1})\circ\mathcal F^{-1} \hat z_j^k)+\alpha_2\mathrm{diag}(M_2^k)\circ u^k}{\beta \sum_j (\mathcal S_j^T | w^{k+1}|^2) +\alpha_2\mathrm{diag}(M_2^k)}; C_u \Big).
\end{aligned}
\right.
\end{equation}
For Step 3, denoting  
\[
z^+=\mathcal A( w^{k+1},u^{k+1})-\tfrac{\Lambda^k}{\beta},
\]
one has 
\[
z^{k+1}=\arg\min_z \tfrac12\langle |z|^2+\varepsilon\mathbf 1_m-(f+\varepsilon\mathbf 1_m)\circ\log(|z|^2+\varepsilon\mathbf 1_m), \mathbf 1_m \rangle+\tfrac{\beta}{2}\|z-z^+\|^2.
\]
The solution can be expressed as
\begin{equation}\label{eqZ-1}
z^{k+1}=\rho^\star \circ\mathrm{sign}(z^+),
\end{equation}
 where $\rho^\star(t)$ was solved by
the gradient projection scheme expressed as:
\begin{equation}\label{eqZ-4}
x_{l+1}=\max\left\{0,x_l-\delta \big((1+\beta-\tfrac{f(t)+\varepsilon}{|x_l|^2+\varepsilon})x_l-\beta z^+(t) \big)\right\}, \forall~l=0,1,\ldots,  \end{equation}
if using the pIPM, or 
\begin{equation}\label{eqZ-3}
x_{l+1}=\max\left\{0, x_l-\delta \big((1+\beta-\tfrac{\sqrt{f(t)+\varepsilon}}{\sqrt{|x_l|^2+\varepsilon}})x_l-\beta z^+(t) \big)\right\}, \forall~l=0,1,\ldots,  
\end{equation}
if using the pAGM with the stepsize $\delta>0,$ and   $x_0:=|z^k(t)|.$
Note that with the penalization parameter $\epsilon=0$, one can directly get the closed form solution by \eqref{eqProxCF} as \cite{wen2012,chang2016Total}.

Under the condition of  sufficient overlapping scan, and bounded preconditioning matrices $M_1^k$ and $M_2^k$, the convergence of the ADMM  can be derived  on the sense that  the iterative sequence generated by above algorithm  converges to a stationary point of the augmented Lagrangian by letting the parameter $\beta$ sufficiently large.

From the point of view of fixed pint analysis, for nonblind problems (knowing the probe $w$), the authors \cite{fannjiang2020fixed} presented a variant ADMM to solve the following optimization problem
\begin{equation}\label{model: ADMM-variant}
\min_{z} \mathcal M(|z|^2;f)+\mathbb I_{X}(z),    
\end{equation}
with $X\subset\mathbb C^m$ defined in \eqref{eq:bilinearSet}.
By introducing the auxiliary variable $\bar z=z$ and decompose the objective functions, the ADMM was proposed in \cite{fannjiang2020fixed} to solve 
\begin{equation}\label{model: ADMM-variant-constr}
\min_{z,\bar z} \mathcal M(|z|^2;f)+\mathbb I_{X}(\bar z), \ s.t. \ z-\bar z=0.    
\end{equation}
To further apply the idea to the BPR, alternating minimization was further adopted as
\begin{equation}\label{alg:ADMM-variant}
\begin{split}
&z^{k+1/2}:=\arg\min_{z} \mathcal M(|z|^2;f)+\mathbb I_{X_1^k}(z);\\
&z^{k+1}:=\arg\min_{z} \mathcal M(|z|^2;f)+\mathbb I_{X_2^k}(z);\\
\end{split}
\end{equation}
where  these two subproblems can be solved via ADMM as inner loop.
Here one has to adjust the constraint sets with the update  probe and sample, i.e. 
\[
\begin{aligned}
&X^k_1:=\{z:~z=\mathcal A(w^k,u)\forall u\in\mathbb C^n\},\\
&X^k_2:=\{z:~z=\mathcal A(w,u^{k+1})\forall w\in\mathbb C^{\bar m}\}.
\end{aligned}
\]
Note that the probe and sample can be readily determined by solving the least squares problem
as
\begin{equation}
\label{eq:LS}
\begin{split}
&u^{k+1}=\arg\min_u \|z^{k+1/2}-\mathcal A(w^k,u)\|^2,\\
&w^{k+1}=\arg\min_w \|z^{k+1}-\mathcal A(w,u^{k+1})\|^2,\\
\end{split}
\end{equation}
which can be solved by \eqref{eq:subsolverAP}.
Although the algorithms worked well with suitable initialization as reported in \cite{fannjiang2020fixed}, the theoretical  convergence  for the blind recovery is still open.

\subsection{Convex programming}
Ahmed et al. \cite{ahmed2018blind} proposed a  convex relaxation  based on a lifted matrix recovery formulation that allows a nontrivial convex relaxation of the Convolution PR.

Consider the Convolution PR as 
\[
f^{Cov}=|\mathcal F \kappa \circ \mathcal F  u|^2.
\] 
One basic assumption for unique recovery is that the variables $\kappa$ and $u$ belong to the subspace of $\mathbb C^n$, i.e. 
\[
\kappa=\mathbf B \bm h,~~u=\mathbf C \bm m,
\]
where $\bm h\in\mathbb C^{k_1}$ and  $\bm m\in\mathbb C^{k_2}$
with known matrices $\mathbf B\in\mathbb C^{n,k_1}$ and $\mathbf C\in\mathbb C^{n,k_2}$ ($k_1, k_2\ll n$).
Then one is concerned with the following problem with  $\bm h,\bm m$ as unknowns
\begin{equation}\label{eq:convex_template}
f^{Cov}=\frac{1}{\sqrt{n}}|\hat {\mathbf B} \bm h\circ \hat{\mathbf C}\bm m|^2
\end{equation}
with $\hat{\mathbf B}:=\sqrt{n}\mathcal F\mathbf B, \hat{\mathbf C}:=\sqrt{n}\mathcal F\mathbf C.$
Further by the lifting technique in semidefinite programming (SDP), the above problem reduces to
\begin{equation}\label{eq:liftMatrix}
f^{Cov}(l)=\frac{1}{n} \langle {\bf b}_l {\bf b}_l^*, \bm H\rangle\langle {\bf c}_l {\bf c}_l^*, \bm M \rangle    
\end{equation}
where 
$\bm H:=\bm h\bm h^*$, $\bm M:=\bm m\bm m^*$ (rank-1 matrices), and
$\langle\cdot,\cdot\rangle$ denotes the Frobenius inner product (trace of multiplication of two matrices).  Here $\bf b_l^*$ and $\bf c_l^*$ are the rows of $\hat {\bf B}$ and $\hat {\bf C}$ respectively. 
By using  a
nuclear-norm minimization, to convexify the rank of matrix, and further transform \eqref{eq:liftMatrix} to a convex constraint, 
then the following convex optimization model  can be derived as
\begin{equation}
\label{model:convex}
\begin{split}
&\min_{\bm H\succcurlyeq \boldsymbol{0}, \bm M\succcurlyeq \boldsymbol{0}}\qquad  {\mathrm Tr}(\bm H)+{\mathrm Tr}(\bm M)\\
&s.t. \qquad\qquad \langle {\bf b}_l {\bf b}_l^*, \bm H\rangle\langle {\bf c}_l {\bf c}_l^*, \bm M \rangle    \geq  \bar f(l),\ \ 0\leq l\leq n-1,
\end{split}
\end{equation}
with $\bar f:=n f^{Cov}$. 

An ADMM  scheme was further developed \cite{ahmed2018blind} to solve \eqref{model:convex}. 
By introducing the convex constraint set
\[
\mathscr C:=\{(\bm v_1, \bm v_2):~\bm v_1(l)\bm v_2(l)\geq \bar f(l), \bm v_1(l)\geq 0~\forall 0\leq l\leq n-1  \}
\]
and $\bm H'=\bm H, \bm M'=\bm M,$
an equivalent form can be given as
\[
\begin{aligned} 
&\min_{\bm H, \bm H', \bm M, \bm M', \bm v_1, \bm v_2}   \mathbb I_{\mathscr C}  (\bm v_1, \bm v_2)+{\mathrm Tr}(\bm H)+{\mathrm Tr}(\bm M)&\\
&\qquad\qquad\qquad\qquad\qquad\qquad+\mathbb I_{\{X\succcurlyeq \boldsymbol{0}\}}(\bm H')+\mathbb I_{\{X\succcurlyeq \boldsymbol{0}\}}(\bm M'),\\
&\quad\quad s.t. \qquad\qquad \bm v_1(l)-\langle {\bf b}_l {\bf b}_l^*, \bm H\rangle=0,\ \  \bm v_2(l)-\langle{\bf c}_l {\bf c}_l^*, \bm M \rangle=0 ,\\
&\qquad\qquad\qquad\qquad \bm H'-\bm H=0,\ \ \  \bm M-\bm M'=0.
\end{aligned}
\]

With the multipliers $\Lambda_k$ for $k=1,2,3,4$ for the totally four constraints,
the augmented Lagrangian with scalar form has the following form
\begin{equation}
\begin{split}
&\mathcal L_c(\bm H, \bm H', \bm M, \bm M', \bm v_1, \bm v_2; \{\Lambda_k\}_{k=1}^4)\\
:=&    \mathbb I_{\mathscr C}  (\bm v_1, \bm v_2)+{\mathrm Tr}(\bm H)+{\mathrm Tr}(\bm M)+\mathbb I_{\{X\succcurlyeq \boldsymbol{0}\}}(\bm H')+\mathbb I_{\{X\succcurlyeq \boldsymbol{0}\}}(\bm M')\\
&+\beta_1\sum_l\left(\langle \Lambda_1(l), \bm v_1(l)-\langle {\bf b}_l {\bf b}_l^*, \bm H\rangle\rangle+\tfrac{1}{2} \|\bm v_1(l)-\langle {\bf b}_l {\bf b}_l^*, \bm H\rangle\|^2\right)\\
&+\beta_1\sum_l\left(\langle \Lambda_2(l), \bm v_2(l)-\langle{\bf c}_l {\bf c}_l^*, \bm M \rangle\rangle+\tfrac{1}{2} \|\bm v_2(l)-\langle{\bf c}_l {\bf c}_l^*, \bm M \rangle\|^2\right)\\
&+\beta_2\langle \Lambda_3,\bm H'-\bm H \rangle+\tfrac{\beta_2}{2}\|\bm H'-\bm H\|^2\\
&+\beta_2\langle \Lambda_4,\bm M'-\bm M \rangle+\tfrac{\beta_2}{2}\|\bm M'-\bm M\|^2,
\end{split}    
\end{equation}
with two positive scalar parameters $\beta_1, \beta_2$.
Then with alternating minimization and update of dual variables $\Lambda_k$, the iterative scheme is obtained. First one can optimize the variable $\bm H$ and $\bm M$ in parallel, and only  consider
\[
\begin{aligned}
\bm H^\star:=&\arg\min_{\bm H} {\mathrm Tr}(\bm H)+\beta_1\sum_l\langle \Lambda_1(l), \bm v_1(l)-\langle {\bf b}_l {\bf b}_l^*, \bm H\rangle\rangle
\\
&+\tfrac{\beta_1}{2} \|\bm v_1(l)-\langle {\bf b}_l {\bf b}_l^*, \bm H\rangle\|^2
+\beta_2\langle \Lambda_3,\bm H'-\bm H \rangle+\tfrac{\beta_2}{2}\|\bm H'-\bm H\|^2.
\end{aligned}
\]
By considering the first order optimality condition  (taking the derivative of the objective function w.r.t. $\bm H$), one obtains 
\[
\mathrm{vec}(\bm H^\star)={\bf T}_1^{-1}\mathrm{vec}\big(\beta_1\sum_l (\bm v_1(l)+\Lambda_1(l)) {\bf b}_l {\bf b}_l^*+\beta_2 (\bm H'-\Lambda_3)-\mathbf I\big),
\]
with 
\[
{\bf T}_1:=\beta_1 \sum_l \mathrm{vec}({\bf b}_l {\bf b}_l^*)\mathrm{vec}({\bf b}_l {\bf b}_l^*)^*+\beta_2\mathbf I.
\]
Similarly, one can determine the optimal  $\bm M^\star$ for the subproblem w.r.t. $\bm M$ by
\[
\mathrm{vec}(\bm M^\star)={\bf T}_2^{-1}\mathrm{vec}\big(\beta_1\sum_l (\bm v_2(l)+\Lambda_2(l)) {\bf c}_l {\bf c}_l^*+\beta_2 (\bm M'-\Lambda_4)-\mathbf I\big),
\]
with 
\[
{\bf T}_2:=\beta_1 \sum_l \mathrm{vec}({\bf c}_l {\bf c}_l^*)\mathrm{vec}({\bf c}_l {\bf c}_l^*)^*+\beta_2\mathbf I.
\]

For the $\bm H'-$subproblem, denoting $\tilde {\bm H}:=\bm H-\Lambda_3$,  one considers the problem
\begin{equation}\label{sub:H}
\bm H'^\star:=\arg\min_{\bm H'} \mathbb I_{\{X\succcurlyeq \boldsymbol{0}\}}(\bm H')+\tfrac{1}{2}\|\bm H'-\tilde {\bm H}\|^2. 
\end{equation}
with the Hermitian matrix $\tilde {\bm H}$(If initializing the multipliers $\Lambda_3$ and $\Lambda_4$ with Hermitian matrices, it can be readily guaranteed that all iterative sequences of these two multipliers are Hermitian).
The closed form solution of \eqref{sub:H} can be directly given as 
\[
\bm H'^\star=\mathrm{Proj}_+(\tilde{\bm H}),
\] 
with the
operator $\mathrm{Proj}_+$ defined as as 
\[
\mathrm{Proj}_+(\tilde{\bm H}):= \bm U \mathrm{diag}(\max\{\mathrm{diag}(\bm \Sigma),0\}) \bm U^*
\] 
and
 $\tilde {\bm H}$ has the eigen-decomposition as  
$\tilde {\bm H}= \bm U \bm \Sigma \bm U^*$
with unitary matrix $\bm U$ and diagonal matrix $\bm\Sigma$.

Similarly 
\begin{equation}\label{sub:M}
\bm M'^\star:=\arg\min_{\bm M'} \mathbb I_{\{X\succcurlyeq \boldsymbol{0}\}}(\bm M')+\tfrac{1}{2}\|\bm M'-(\bm M-\Lambda_3)\|^2.
\end{equation}
One can directly get the closed form solution
\[
\bm M'^\star=\mathrm{Proj}_+(\bm M-\Lambda_3).
\]

The subproblems w.r.t the variables $\bm v_1,$ and  $\bm v_2$ can be solved in an element-wise manner, due to the independence of the optimization problem for each element of these two variables. 
Since they can be derived with standard discussion based on Karush-Kuhn-Tucker optimality conditions,   the details here are omitted. 

Hence  all procedures to get the iterative scheme  are summarized by further combined with the update of the multipliers 
as 
\[
\begin{split}
&\Lambda_1(l)\leftarrow \Lambda_1(l) +\bm v_1(l)-\langle {\bf b}_l {\bf b}_l^*, \bm H\rangle;\\
&\Lambda_2(l)\leftarrow \Lambda_2(l) +\bm v_2(l)-\langle {\bf c}_l {\bf c}_l^*, \bm M\rangle;\\
&\Lambda_3\leftarrow\Lambda_3+\bm H'-\bm H;\\
&\Lambda_4\leftarrow\Lambda_4+\bm M'-\bm M.
\end{split}
\]
Please see more details in the appendix of \cite{ahmed2018blind}.

As reported in \cite{ahmed2018blind}, this convex method showed excellent agreement with the theorem in the case of random subspaces.
However, it was less effective on  deterministic subspaces, including
partial Discrete Cosine Transforms or partial Discrete Wavelet Transforms. One should also notice that although the model is convex, the lifting technique increased  the dimension of original nonconvex optimization problem greatly, at the order of square of the original dimension, 
causing huge memory requirement as well as computational complexity. That may limit practical applications, especially for reconstructing 2D  images or volumes.

It seems rather difficult to adopt the same convex method for other cases of BPR, since they cannot be rewritten as the same form as \eqref{eq:convex_template}. Convexifying a general BPR problem  should be an interesting research direction in future.

\subsection{Second order algorithm using Hessian}

The second order algorithms  relying on the Hessian of the nonlinear optimization problems,  have also been developed for PR problem, such as using Newton method (NT) \cite{qian2014efficient,yeh2015experimental}, Levenberg-Marquardt method (LM) \cite{ma2018globally,kandel2021efficient} or Gauss-Newton algorithm (GN)\cite{gao2017phaseless}. 
Consider the following problem by rewriting \eqref{eqPB}
\begin{equation}\label{eq:nonblindLS}
    \min_u \mathcal Q(u),
\end{equation}
with $\mathcal Q(u):=\mathcal M(| A_w u|^2,f)$.
Given the initial guess $u^0$, 
\begin{equation}\label{eq:Newton}
\mathcal Q(u)\approx\mathcal Q(u^0)+\Re\langle \nabla_u \mathcal Q(u^0), u-u^0 \rangle+\frac{1}{2}\Re(\langle \nabla^2_u \mathcal Q(u^0) (u-u^0), u-u^0 \rangle),
\end{equation}
where $\nabla_u^2$ denotes the Hessian operator.
Then a new estimate $u^1$ for the stationary point can be obtained by solving the following systems
\[
\nabla^2_u \mathcal Q(u_0) (u^1-u^0)=-\nabla_u \mathcal Q(u^0).
\]
Assuming the Hessian matrix is nonsingular, the iterative scheme by NT  is derived as
\begin{equation}
\text{Newton method:}\qquad u^{k+1}=(\nabla^2_u \mathcal Q(u^k))^{-1}\big( u^k-\nabla_u \mathcal Q(u^k)\big)~\forall k.
\end{equation}
The gradient  is given below:
\begin{equation}
    \label{eq:gradient}
\nabla_u \mathcal Q(u)=
\left\{
\begin{aligned}
&A_w^*\big(A_w u-\tfrac{\sqrt{f}}{|A_w u|}\circ A_w u\big); &\text{(AGM)}\\
&A_w^*\big(A_w u-\tfrac{f}{|A_w u|^2}\circ A_w u \big); &\text{(IPM)}\\
&2A_w^*\big(|A_w u|^2\circ A_w u-f\circ A_w u \big); &\text{(IGM)}
\end{aligned}
\right.
\end{equation}
where the objective function $\mathcal Q(u)$  is rewritten as  $\mathcal Q(u)=\mathcal M(|A_w u|^2, f)$ by denoting the matrix $A_w$ as \eqref{eq:linearOperator} and the detailed forms of the operators can be found in \eqref{eq:AwT}-\eqref{eq:AtAu}.
The Hessian matrices for three metrics are complicated, and please see Appendix A of \cite{yeh2015experimental}.

 More efficient algorithms including LM and GN were developed,  concerned with the nonlinear least square problems (NLS) \eqref{eq:nonblindLS} with the AGM and IPM metrics (Please see \eqref{eqDF}.
Namely, by denoting the residual function
\[
\bm r(u)=\left\{
\begin{aligned}
&|A_w u|-\sqrt{f}; &\text{(AGM)}\\
&|A_w u|^2-f;      &\text{(IGM)}\\
\end{aligned}
\right.
\]
consider the NLS problem below:
\[
\min_u \mathcal Q(u)=\frac12\|\bm r(u)\|^2.
\]
Then with Jacobian matrix  as
\[
\bm J(u):=\nabla_u \bm r(u)=\left\{
\begin{aligned}
&\mathrm{diag}(\mathrm{sign}(\mathrm{conj}(A_w u)))A_w; &\text{(AGM)}\\
&\mathrm{diag}(\mathrm{conj}(A_w u))A_w &\text{(IGM)};\\
\end{aligned}
\right.
\]
the GN method considered  
\[
GN(u):=\bm J^*(u)\bm J(u),
\]
as an estimate of the Hessian matrix, that leads to the following iterative scheme
\begin{equation}
\begin{aligned}
\text{Gauss-Newton method:}\qquad u^{k+1}&=(\bm J^*(u^k)\bm J(u^k))^{-1}\big( u^k-\nabla_u \mathcal Q(u^k)\big)\\
                                         &=(\bm J^*(u^k)\bm J(u^k))^{-1}\big( u^k-\bm J^*(u^k)\bm r(u^k)\big)~~~\forall k.
\end{aligned}
\end{equation}
Gao and Xu \cite{gao2017phaseless} further proposed a global convergent GN algorithm with re-sampling for PR problem, which partial phaseless data was used to reformulate the GN matrix and the gradient per loop. 

The Hessian matrix or the GN matrix  can not guaranteed to be nonsingular  practically. Hence the LM method interpreted as a regularized variant of GN was proposed as
\begin{equation}
\begin{aligned}
\text{LM  method:}\qquad u^{k+1}&=(\bm J^*(u^k)\bm J(u^k)+\mu^k\mathbf I)^{-1}\big( u^k-\bm J^*(u^k)\bm r(u^k)\big)~~~\forall k.
\end{aligned}                                         
\end{equation}
with the adaptive parameter $\mu^k$. Readily one knows $\mu^k$ cannot be too large, otherwise the Hessian information is useless.  To obtain fast convergence,
Marquardt \cite{Marquardt1963} proposed the following strategy for $\mu^k$ depending on the diagonal matrix of GN matrix as
\[
\mu^k=\mu_0 \bm {D_g}(\bm J^*(u^k)\bm J(u^k))),
\]
with $\bm D_g(A)$ denoting the  diagonal matrix with the elements from the main diagonal of the matrix $A$.
Yamashita \cite{Yamashita2001}, and Fan$\&$Yuan  \cite{fan2005quadratic} proposed the scheme depending on the objective function value below:
\begin{equation}\label{eq:FY}
\mu^k=(\mathcal Q(u^k))^{\frac{\nu}{2}}.
\end{equation}
with $\nu\in[1,2].$
 Ma, Liu and Wen \cite{ma2018globally}  further improved the scheme \eqref{eq:FY} as choosing a larger value   
when the iterative solution $u^k$ is far away from the global minimizer, i.e.
\[
\mu^k=\mathrm{Thresh}(u^k)(\mathcal Q(u^k))^{\frac{\nu}{2}},
\]
with 
\[
\mathrm{Thresh}(u)=\left\{
\begin{aligned}
&\tau, &\text{if~}\mathcal Q(u^k)\geq c_0 \|u^k\|^2,\\
& 1, &\text{otherwise,}
\end{aligned}
\right.
\]
 with $\tau\gg 1$ and parameter $c_0>0.$

The mentioned algorithms including  \cite{qian2014efficient,gao2017phaseless,ma2018globally}  focused on nonblind PR. 
With the generalized GN method and automatic differentiation, Kandel et al. \cite{kandel2021efficient}
proposed a variant LM algorithm for blind recovery, where especially for IPM based metric, it employed the generalized GN (GGN) as 
\[
GGN(u):=\bm J^*(u^k)\nabla^2_g \mathcal M(| A_{u^k} w|^2,f)\bm J(u^k)
\]
with $\mathcal M(g,f)$ defined in \eqref{eqDF}.
Following the same manner with alternating minimization, one can easily derive the second-order algorithm for the blind problem as
\begin{equation}
\begin{aligned}
&u^{k+1}=\arg\min_u \mathcal M(|A_{w^k}u |^2,f);\\
&w^{k+1}=\arg\min_w \mathcal M(| A_{u^{k+1}} u|^2,f);\\
\end{aligned}    
\end{equation}
where the both two subproblems are solved by NT,  GN  or LM algorithms.

\subsection{Subspace method}

The subspace method \cite{sadd2003iterative,CSIAM-AM-2-585} is a very powerful algorithm, iteratively refining the variable in the subspace of solution, which includes the Krylov subspace method as well-known conjugate gradients method, domain decomposition method, and multigrid method.
 It originally focusing on solving the linear equations or least-square problems, now  has been successfully extended to nonlinear equations or nonlinear optimization problems.
 In this part,   the subspace methods for the PR and BPR problems will be reviewed.

\subsubsection{Nonlinear conjugate gradient algorithm}

Consider the following optimization problem
\[
\min f(\bm x).
\]
By the nonlinear conjugate gradient (NLCG) algorithm, 
the iterative scheme can be given below:
\begin{equation}\label{NLCG}
\begin{aligned}
&\bm x^{k+1}=\bm x^k+\alpha^k \bm d^k;\\
&\bm d^k=-\nabla_{\bm x} f(\bm x^k)+\beta^{k-1}\bm d^{k-1}~~\forall  k\geq 1, 
\end{aligned}    
\end{equation}
with the stepsize $\alpha^k$ and  weight $\beta^{k-1}$, 
where $\bm d^k$ is the search direction.
One may notice that the search direction $\bm d^k$ in NLCG is the combination of the  gradient and the search direction $\bm d^{k-1}$ with the weight $\beta^{k-1}.$
To get optimal parameters, the stepsize $\alpha^k$ is selected by the monotone line search procedures, while the weight $\beta^k$ is determined  based on the gradient $\nabla_{\bm x} f(\bm x^{k-1})$, $\nabla_{\bm x} f(\bm x^k)$ and the search direction $\bm d^{k-1}$ (typically five different formulas \cite{CSIAM-AM-2-585}).  

The NLCG has been successfully applied to the BPR problem \cite{thibault2012maximum,qian2014efficient}.
For example, Thibault and Guizar-Sicairos \cite{thibault2012maximum} adopted the NLCG to solve the CDI problem. 
The iterative scheme can be given below:
\begin{equation}\label{NLCG-ptycho}
\begin{aligned}
&(w^{k+1}, u^{k+1})=(w^{k}, u^{k})+\alpha^k \bm \Delta^k;\\
&\bm\Delta^k=-\bm g^k+\beta^{k-1}\bm\Delta^{k-1}~~\forall  k\geq 1, 
\end{aligned}    
\end{equation}
with the gradient $\bm g^k:=(\nabla_w \mathcal M(|\mathcal A(w^k,u^k)|^2,f),\nabla_u \mathcal M(|\mathcal A(w^k,u^k)|^2,f))$ calculated by \eqref{eq:gradient} and $\bm \Delta:=(\bm \Delta_w, \bm \Delta_u)$.
The weight $\beta^{k-1}$ is derived by the Polak-Rib\'ere formula as 
\[
\beta^{k-1}=\frac{\langle \bm g^k, \bm g^k \rangle-\Re(\langle \bm g^k, \bm g^{k-1} \rangle) }{\langle \bm g^{k-1}, \bm g^{k-1} \rangle}.
\]
To further get $\alpha^k$,  by estimating  
$\mathcal M(|\mathcal A(w^k+\alpha \bm \Delta_w,u^k+\alpha \bm \Delta_u)|^2,f)$ by the low-order polynomial as
\[
 \mathcal M(|\mathcal A(w^k+\alpha \bm \Delta_w,u^k+\alpha \bm \Delta_u)|^2,f)\approx \sum_{t=0}^8 c_t\alpha^t,
\]
 the authors adopted the Newton–Raphson algorithm in order to minimize the following one-dimension problem:
\[
\alpha^k:=\arg\min_\alpha \sum_{t=0}^8 c_t\alpha^t.
\]

\subsubsection{Domain decomposition method}
The domain decomposition methods (DDMs) allow for highly parallel computing with good load balance,  by decomposing the equations on whole domain to the problems on relatively small subdomains with information synchronization on the partition interfaces. They have played a great role in solving partial differential equations numerically  and recently  been successfully extended to large-scale image restoration, image reconstruction and other inverse problems, e.g. \cite{xu2010two,Chang2015sims,langer2019overlapping,lee2019finite,chang2021over}  and references therein. 
For ptychography imaging, several parallel algorithms \cite{nashed2014parallel,guizar2014high,marchesini2016sharp,enfedaque2019high,chang2021over} have been developed.
 Specially, for convention ptychography, Chang et al. \cite{chang2021over} proposed an overlapping DDM with the ST-AGM as defined in \eqref{eq:STAGM}, with fewer communication cost and theoretical convergence guarantee.

First give the domain decomposition. Denote the  whole region $\Omega:=\{0, 1, 2, \cdots, n-1\}$ in the discrete setting. There exists the two-subdomain overlapping DD $\{\Omega_d\}_{d=1}^2$,  such that  \[
\Omega=\bigcup_{d=1}^{2} \Omega_d\]
with $\Omega_d:=\{l^d_0, l^d_1,\cdots, l^d_{n_d-1}\}$, and the overlapping region is denoted as \[
\Omega_{1,2}:=\Omega_1\cap\Omega_2=\{\hat l_0, \hat l_1, \cdots, \hat l_{\hat n-1}\}.\]
Here consider a special overlapping DD as shown in Fig. \ref{fig2}. Denote the restriction operators $R_1, R_2$ as 
\[
R_d u=u|_{\Omega_d},\ \ \  R_{1,2}u=u|_{\Omega_{1,2}},\]
i.e. 
\[
\begin{aligned}
&(R_d u)(j)=u(l^d_j)\forall \ 0\leq j\leq n_d-1,\\
&(R_{1,2}u)(j)=u(\hat l_j)\forall 0\leq j\leq \hat n-1.
\end{aligned}
\]
Then two groups of localized shift operators can be introduced $\{\mathcal S_{j_d}^d\}_{j_d=0}^{J_d-1}, d=1,2$, with $\sum_d J_d=J.$

For nonblind problem,   denote the linear operators $A_1, A_2$  on the subdomains as
\begin{equation}
A_d u_d:=\big((\mathcal F(w\circ\mathcal S_0^d u_d))^T, (\mathcal F(w\circ\mathcal S_1^du_d))\big)^T,\cdots, (\mathcal F(w\circ\mathcal S_{J_d-1}^du_d))^T)^T\ \ d=1,2.
\label{eq:Ad_def}
\end{equation}
Naturally, the measurement $f$ is also decomposed to two non-overlapping parts $f_1, f_2$, i.e. 
\[
f_d:=|A_d R_d u|^2
\]
Hereafter, consider the nonlinear optimization problem with ST-AGM. In order to enable the parallel computing of $u_1$ and $u_2$,  introduce an auxiliary variable $v$ which is only defined {in the overlapping region $\Omega_{1,2}$,} and then are concerned with the following model:
\begin{equation}
\begin{aligned}
    &\min\nolimits_{u_1, u_2, v}&\sum\nolimits_{d=1}^2\mathcal G_\epsilon(A_d u_d;f_d)\\ &\text{~~~~s.t.}& \qquad\mathcal {\mathbf  \pi}_{d,3-d}u_d-v=0,~ d=1,2.
\end{aligned}
\label{model:twosubdomain}
\end{equation}
In order to develop an iterative scheme without inner loop as well as with fast convergence for large-step scan,} two  auxiliary variables $z_1, z_2$ are introduced below:
\begin{equation}
\begin{aligned}
    &\min\nolimits_{u_1, u_2, v, z_1, z_2}&&\sum\nolimits_{d=1}^2\mathcal G_\epsilon(z_d;f_d)\ \\
& s.t. &&  \mathcal {\mathbf  \pi}_{d,3-d}u_d-v=0, \ \ \ A_d u_d-z_d=0,\ d=1, 2.
\end{aligned}
\label{eq:model-constraint}
\end{equation}
Then it is quite standard to solve the saddle point problem by ADMM. The details are omitted here, and please see more details in \cite{chang2021over}.

\begin{figure}[ht]
\begin{center}
\subfigure[]{\includegraphics[width=.4\textwidth]{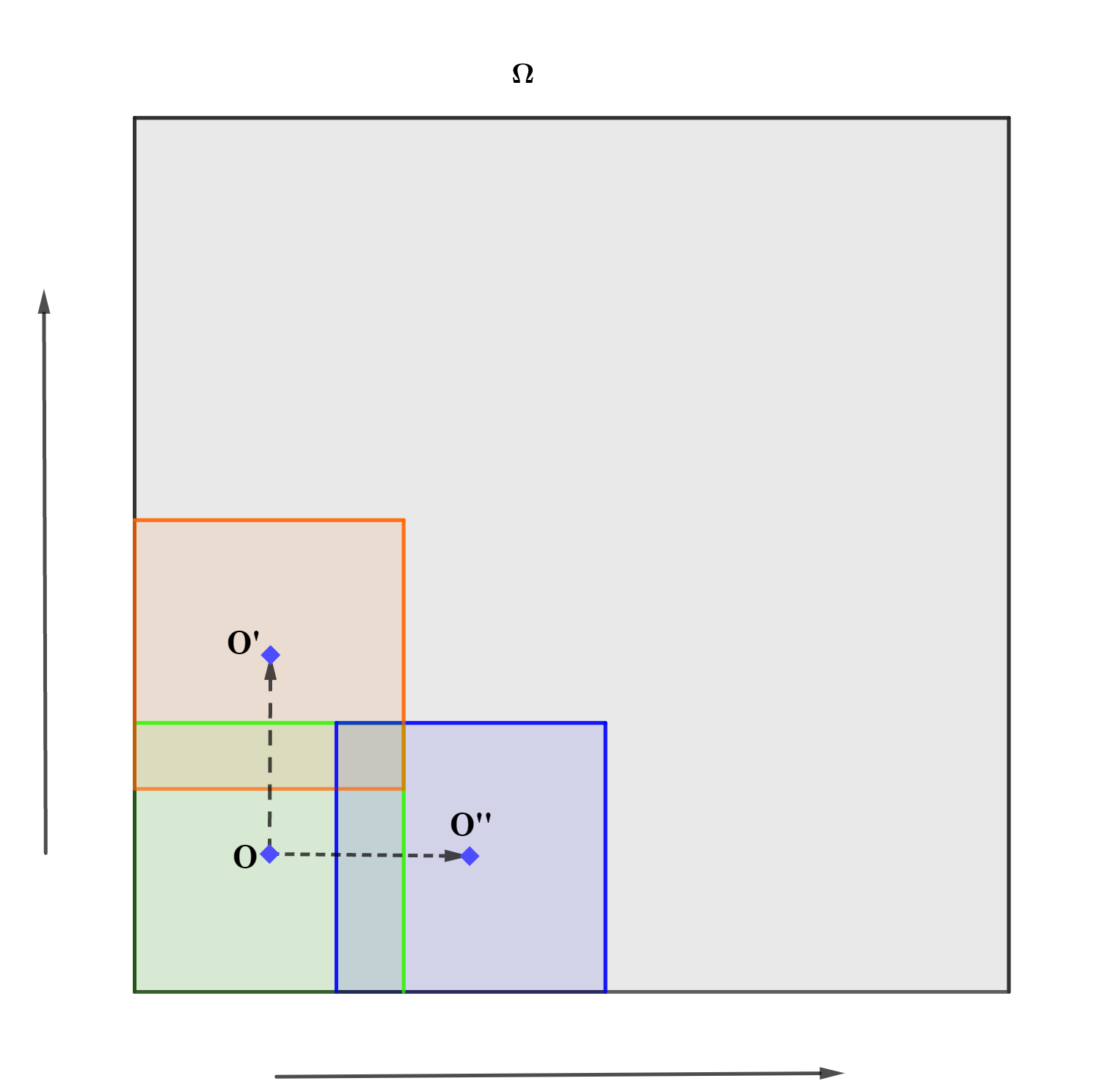}}
\subfigure[]{\includegraphics[width=.4\textwidth]{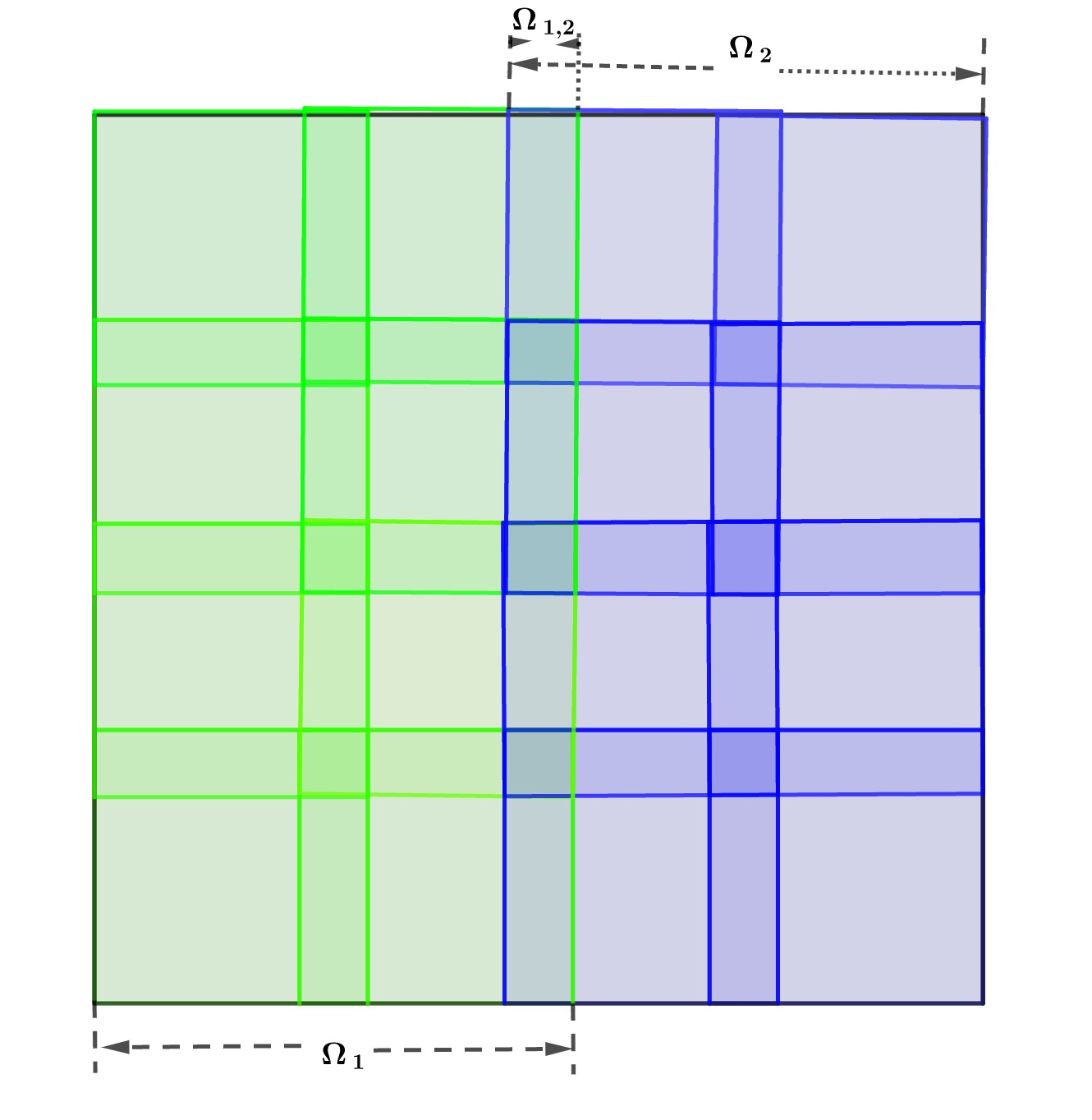}}
\end{center}
    \caption{(a) Ptychography scan in the domain $\Omega$ (grid scan): the starting scan centers at point {\bf O}, and then move up (or to the right) with the center point {\bf O'} (or {\bf O"} ); (b) Two-subdomain DD (totally $4\times 4$ frames): The subdomains $\Omega_1, \Omega_2$ are generated by two $4\times 2$-scans, and the overlapping region $\Omega_{1,2}=\Omega_1\cap\Omega_2$. }
\label{fig2}
\end{figure}

Then for blind recovery,  in order to reduce the grid pathology \cite{chang2018Blind} (ambiguity derived by the multiplication of any periodical function and the true solution) due to grid scan, introduce the support set constraint of the probe, i.e. $\mathcal O:=\{w:\ (\mathcal Fw)(j)=0,\ j\in \mathscr J\},$ { with the support set $\bar {\mathscr J}$ denoted as the complement of the set $\mathscr J$ (index set for zero values for the Fourier transform of  the probe). } 
Then  consider the blind ptychography problem for two-subdomain DD:
\begin{equation*}
\begin{aligned}
    &\min\nolimits_{\{w, u_1, u_2, v\}}\sum\nolimits_{d=1}^2\mathcal G_\epsilon(\mathcal A_d(w, u_d);f_d)+\mathbb I_{\mathcal O}(w),\\ 
    &s.t.\qquad {\pi}_{d,3-d}u_d-v=0, d=1,2, \\
\end{aligned}
\end{equation*}
where the bilinear mapping
$\mathcal A_d(w, u_d)$ is denoted as
\[(\mathcal A_d)_{j_d}(w,u_d):=\mathcal F (w\circ \big(S^d_{j_d}u_d))\ \forall\ 0\leq j_d\leq J_{d}-1,
\]
with
$\sum_{d=1}^{2} J_{d}=J$, and the indicator function $\mathbb I_{\mathcal O}$. 
To enable parallel computing, consider the following constraint optimization problems
\begin{equation*}
\begin{aligned}
    \min_{\{w, w_1, w_2, u_1, u_2, v, z_1, z_2\}}&\sum_{d=1}^2\mathcal G_\epsilon(z_d;f_d)+\mathbb I_{\mathcal O}(w)\ \\
s.t.\qquad\qquad& \ \mathcal {\pi}_{1,2}u_1-v=0,\ \  \mathcal { \pi}_{2,1}u_2-v=0, \\
& w_d=w, \ z_d=\mathcal A_d(w_d, u_d), d=1,2,  \ \
\end{aligned}
\end{equation*}
which was also efficiently solved by ADMM.

\subsubsection{Multigrid methods}
The multigrid method (MG) is a standard framework in order to accelerate solving partial differential equations \cite{hackbusch2013multi},  large-scale linear equations \cite{xu_zikatanov_2017} and related optimization problems  \cite{borzi2009multigrid} with the full approximation scheme (FAS) \cite{brandt2011multigrid}. 
A multigrid-based optimization framework based on \cite{nash2000multigrid} to reduce the computational
for nonblind ptychographic phase retrieval was proposed by Fung and Di \cite{fung2020multigrid}, which utilized the
hierarchical structures of the measured data.

Consider the following feasible problem \cite{fung2020multigrid} as
\begin{equation}
\min_u \sum_j \| \mathcal F^*(\sqrt{f_j}\circ \mathrm{sign}(\mathcal F(w\circ \mathcal S_j u)))-w\circ \mathcal S_j u \|^2,     
\end{equation}
which is  equivalent to the problem
\[
\min_u \mathcal M(|\mathcal A(w,u)|^2,f),     
\]
with AGM metric.
Then the multigrid optimization framework based on FAS was further developed,  where the coarse-grid subproblem was interpreted as a first-order approximation to the fine grid problem. However, it is unclear that how to extend the current algorithm to the blind problem.

\section{Discussions}\label{sec4}

\subsection{Experimental  issues}

\subsubsection{Probe drift}
Probe drift happens in ptychography, when the data is very noisy. The mass center of the iterative probe will eventually touch the boundary such that the iterative algorithms fail eventually. Hence the joint reconstruction will cause instability of the iterative algorithms from noisy experimental data. One simple strategy proposed by \cite{marchesini2016sharp} is to shift the probe to the mass center of the complex image periodically. 
Other possible way is to assume that  the compact support condition for the probe, or  to get additional measurement for the probe by letting the light go through the vacuum as \cite{marchesini2016sharp,chang2018Blind}. The related numerical stability shall be investigated, and one can refer to Refs. \cite{huang2020estimation,huang2021uniqueness} for nonblind PR.    

\subsubsection{Flat samples}
When the sample is nearly flat (such as weak absorption or scattering for biological specimens using  hard x-ray sources), there will be no sufficient diversity of the  measured phaseless data even by very dense scan. In such case, the iterative algorithms  mentioned in this survey will become 
slow and the recovered image quality get worse. Acquiring of scattering map by linearization for  large features of the sample \cite{Dierolf_2010ptycho} or modeling with additional Kramers-Kronig relation (KKR) \cite{Hirose:17}, were exploited to improve the reconstruction quality.   Besides, pairwise relations between adjacent frames were considered in \cite{marchesini2014rank} to accelerate  projection algorithms for the flat sample.

\subsubsection{Background retrieval}
 Parasitic scattering termed as background often happens experimentally, which may come from any element along the beam path  other than the sample and the optical elements desired harmonic order \cite{chang2019advanced}. Direct reconstruction without background removal will introduce structural artifacts to the reconstruction images. Several methods were designed, such as 
preconditioned gradient descent \cite{marchesini2013augmented}, preprocessing method \cite{wang2017background}, and ADMM for nonlinear optimization method with framewise-invariant background \cite{chang2019advanced}. It is still a challenging problem since the practical background is sophisticated, and cannot be assumed to be framewisely invariant.

\subsubsection{High dimensional problems}
The formula for all four cases for BPR  holds for  a thin (2D) object in paraxial approximation.
For thick samples, the linear propagation as \eqref{equations:blindPR} will cause obvious errors, and one has to consider the nonlinear transform as \cite{dierolf2010ptychographic}. 
Other than the 3D imaging, high dimensional problems may result from the spectromicroscopy \cite{maiden2013soft}, multi-modes decomposition of partial coherence \cite{Thibault2013reconstructing,chang2018partially}, and dichroic ptychography \cite{chang2020analyzer,Lo2021xray}. Such strong nonlinearity coupling with the high dimensional optimization causes difficulties for designing the stable and high-throughput algorithm.

\subsection{Theoretical  analysis}

\subsubsection{Convergence of iterative algorithms }

Other than the projection onto  nonconvex modulus constraint for nonblind PR,  APs \cite{thibault2009probe,marchesini2016sharp} for BPR involve   the bilinear constraint set.
Some progress has been made for the general PR problem using projection algorithms~ \cite{hesse2013nonconvex,marchesini2015alternating,chen2016fourier}. however, the corresponding convergence theories for  BPR are still unclear,  Moreover, only the  PHeBIE \cite{hesse2015proximal} and ADMM based algorithm \cite{chang2018Blind}  for BPR provided rigorous convergence analysis. Hence, it is of great importance to either study the  convergence of existing algorithms, or develop new algorithms with clear convergence guarantee in future.

\subsubsection{Uniqueness analysis}
Uniqueness can be guaranteed for 1-D nonblind pytchographic  PR for nonvanishing signals with the probe of proper size\cite{jaganathan2016stft}. It can also be guaranteed for BPR \cite{bendory2019blind}. By letting two signals lie in low-dimensional random subspaces,  the uniqueness was obtained \cite{ahmed2018blind} with sufficient measurements. For 2D imaging problems, with a randomly
phased probe, the uniqueness can be proved for the measurements which is strongly connected and possess
an anchor.  See more discussions on more general cases together with sparse signals in \cite{Grohs2020phase}.  
Readily for ptychography, nontrivial ambiguity including  periodical function and linear phase  exists for raster scan. Rigorous analysis about more general ambiguity was given \cite{Fannjiang2019raster}.  Experimentally more flexible spiral or random scan \cite{huang2014optimization} have been exploited for stable recovery.

\section{Conclusions}\label{sec5}

In this survey,  a short review of the iterative algorithms is provided for the nonlinear optimization problem arising from the BPR problem, mainly consisting of three types of algorithms as the first-order operator-splitting algorithms, and second order algorithms and subspace methods. There still exist sophisticated experimental issues and challenging  theoretical analysis, which  is further discussed in the last part.   
Learning based methods have been a powerful tool for solving inverse problem and PR problems, which  are not included in this survey.

This survey focuses the BPR problems with forms expressed as \eqref{eq:bilinear}. However,  not all the BPR problem belongs to the categories of \eqref{eq:bilinear}. Very recently, a resolution-enhanced parallel coded ptychography (CP) technique \cite{jiang2021resolution,jiang2022ptychographic} was reported which 
achieves the highest numerical aperture. With  the sample $u$ and the the transmission profile of the engineered surface $w$, the phaseless data was generated as 
\[
f^{CP}_j=\big|\big(w\circ (\mathcal S_j u\circledast\kappa_1)\big)\circledast\kappa_2 \big|^2,
\]
with $\kappa_1, \kappa_2$ as two known PSFs .
Such advanced cases should be further investigated. 

\section*{Acknowledgments}
The work of the first author was partially supported by the NSFC (Nos. 12271404, 11871372, 11501413) and Natural Science Foundation of Tianjin (18JCYBJC16600).
The authors would like to thank Prof. Guoan Zheng  for the helpful discussions.

\bibliographystyle{siamplain}

\bibliography{rD}

\end{document}